# MASS MATRIX INTEGRATION SCHEME FOR FIFTEEN-NODE WEDGE ELEMENT.

E Hanukah

Faculty of Mechanical Engineering, Technion – Israel Institute of Technology, Haifa 32000, Israel
Corresponding author Email: eliezerh@tx.technion.ac.il

**Abstract**

At present, mass matrix of solid fifteen node wedge element is computed by means of eighteen-point (Gauss points) numerical integration scheme. Herein, this widely accepted scheme is being challenged. We derive a novel, easy-to-implement, ten-point integration rule. To this end, the metric (Jacobian determinant) is approximated using special second order interpolation, requiring ten evaluation points. This polynomial approximation permits further analytical integration, which is accompanied by convenient coefficient matrix definition. Coefficient matrices (equivalent to weights), allow the new rule to be formulated in a well-known manner. Preliminary numerical study considering both fine and a coarse mesh is conducted. In fact, significant accuracy superiority over eighteen-point scheme is established for all the coarseness range. In conclusion, our ten-point mass matrix scheme over-performs the standard 18-point integration rule in both the accuracy and in computational effort.

1. **Introduction**

Mass coefficients, internal forces, stiffness matrix, all require integration in the element domain, which is most commonly obtained with the help of numerical integration schemes [1]. Eighteen-point scheme is needed to accurately calculate the consistent mass matrix of fifteen-node wedge element [2]. Several studies exploit the idea of closed-form integration for stiffness matrixes [3-10], significant time savings has been established. Furthermore, hierarchical semi-analytical displacement based approach is used to model three dimensional finite bodies e.g. [11-13] yielding new analytical solutions.

In present study we follow the basic guidelines presented in [14-16]. The metric is approximated using interpolations of different order. Zero order approximation requires metric evaluation at one point, usually the centroid; first order interpolation involves four points and linear interpolation function; quadratic interpolation uses ten evaluation points and second order shape functions, etc. Hence, metric takes simple polynomial series form, which allows further analytical integration. Though, generally speaking, metric is not limited to polynomial interpolations. Coefficient matrices definition allows familiar representation of the resulting integration rule.

We've derived quadratic metric (QM) polynomial interpolation, which together with analytical integration resulted in ten-point integration rule. Preliminary numerical study including fine and coarse mesh elements, revealed significant accuracy superiority over the 18-point scheme. According to our findings, it is beneficial to use our QM mass matrix integration rule over the standard 18-point scheme in terms of accuracy and computational effort.

The outline of the paper as follows. Section 2 recalls important details of the fifteen-node solid wedge element concluding with mass matrix formulation and basic guidelines of the method being used. Natural coordinates, shape functions, metric and jacobian matrix, numerical and analytical integration in the element domain are recalled. Section 3 illustrates the functional form of the jacobian matrix and the metric. Special CM, LM, and QM metric interpolations are proposed. Section 4 presents the resulting integration schemes; coefficient matrices are defined. Section 5 contains preliminary numerical accuracy study stressing out the comparison between the standard 18-point scheme and our ten-point rule. Both, fine and a coarse mesh has been considered. Section 6 summarizes and records our conclusions.

## 2. Background

Initial nodal locations (e.g. [17] pp.75) are denoted by $\mathbf{X}_i$ $(i=1,..,15)$. While $X_{mi}$ $(m=1,2,3, i=1,..,15)$ stand for nodal components in terms of global Cartesian coordinates system $\mathbf{X}_i = X_{mi}\mathbf{e}_m$ $(m=1,2,3, i=1,..,15)$, summation convention on repeated indexes is implied. Here and throughout the text, bold symbols traditionally denote vector or tensor quantities. Natural coordinates $\{\xi,\eta,\zeta\}$ $0\leq\xi,\eta,\zeta\leq 1$ define the domain $0\leq\xi\leq 1-\eta$, $0\leq\eta\leq 1$, $0\leq\zeta\leq 1$ (see Fig. 1)

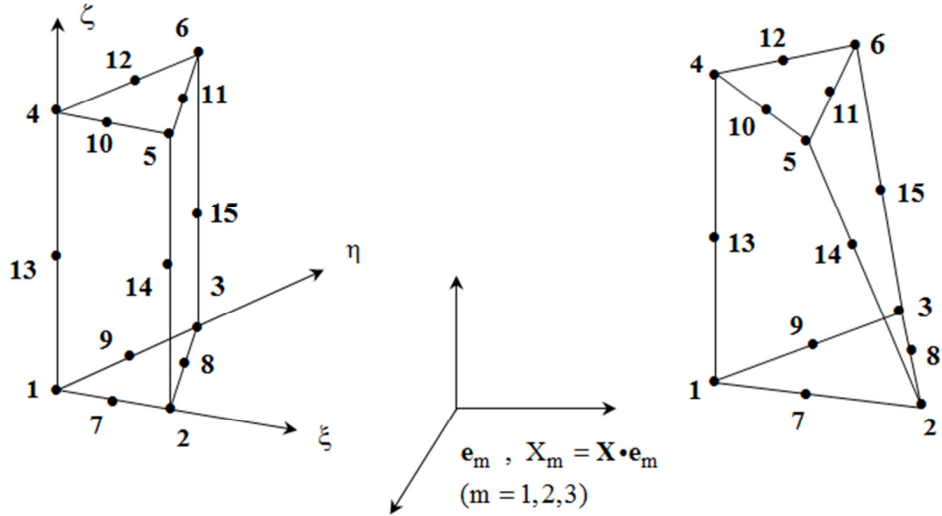

Figure 1: Showing the Global Cartesian coordinate system $\mathbf{e}_m$ $(m=1,2,3)$, natural coordinates $\{\xi,\eta,\zeta\}$, nodes numbering, straight-sided (constant metric/parent) element on the left and general (coarse) mesh on the right.

Shape functions $\varphi^i(\xi,\eta,\zeta)$, $(i=1,..,15)$ are recalled

$$\varphi^1 = -(1-\xi-\eta)(1-\zeta)(2\xi+2\eta+\zeta)/2 \;,\quad \varphi^2 = \xi(1-\zeta)(2\xi-\zeta-2)/2$$
$$\varphi^3 = \eta(1-\zeta)(2\eta-\zeta-2)/2 \;,\quad \varphi^4 = -(1-\xi-\eta)(1+\zeta)(2\xi+2\eta-\zeta)/2$$
$$\varphi^5 = \xi(1+\zeta)(2\xi+\zeta-2)/2 \;,\quad \varphi^6 = \eta(1+\zeta)(2\eta+\zeta-2)/2$$
$$\varphi^7 = 2\xi(1-\xi-\eta)(1-\zeta) \;,\; \varphi^8 = 2\xi\eta(1-\zeta) \;,\; \varphi^9 = 2\eta(1-\xi-\eta)(1-\zeta)$$
$$\varphi^{10} = 2\xi(1-\xi-\eta)(1+\zeta) \;,\; \varphi^{11} = 2\xi\eta(1+\zeta) \;,\; \varphi^{12} = 2\eta(1-\xi-\eta)(1+\zeta)$$
$$\varphi^{13} = (1-\xi-\eta)(1-\zeta^2) \;,\; \varphi^{14} = \xi(1-\zeta^2) \;,\; \varphi^{15} = \eta(1-\zeta^2)$$
(1)

Material point X inside an element domain is denoted by **X**

$$\mathbf{X} = \varphi^i(\xi,\eta,\zeta)\mathbf{X}_i \;,\; (i=1,..,15) \tag{2}$$

Metric or determinant of the Jacobian matrix of global-local coordinate's transformation

$$J = \mathbf{X}_{,1} \times \mathbf{X}_{,2} \cdot \mathbf{X}_{,3} = |J_{mn}| = J_{11}J_{22}J_{33} - J_{11}J_{23}J_{32} - J_{31}J_{22}J_{13} - J_{21}J_{12}J_{33} + J_{21}J_{32}J_{13} + J_{31}J_{12}J_{23}$$
$$J > 0 \;,\; J_{mn}(\xi,\eta,\zeta;X_{ki}) = (\mathbf{X}\cdot\mathbf{e}_m)_{,n} \;,\; (i=1,..,15, m,n,k=1,2,3) \tag{3}$$

Where (×) and (·) stand for vector and scalar products, $|\cdot|$ stand for determinant operator, comma denotes partial differentiation with respect to natural coordinates. Differential volume element is given by

$$dV = J(\xi,\eta,\zeta;X_{mi})d\xi d\eta d\zeta \;,\; (m=1,2,3, i=1,..,15) \;,\; V = \int_V dV \tag{4}$$

Isoparametric formulation (e.g.[18] pp.104) for mass conserving element yield consistent, symmetric, positive definite mass matrix

$$M^{ij} = \int_V \rho_0 \varphi^i \varphi^j dV = \rho_0 \int_0^{+1}\int_0^{+1}\int_0^{1-\eta} \varphi^i \varphi^j J d\xi d\eta d\zeta \;,\; M^{ij} = M^{ij} \;,\; (i,j=1,..,15) \tag{5}$$

Numerical integration of (5) is recalled (e.g. [18] pp.126)

$$\int_V \rho_0 \varphi^i \varphi^j J d\xi d\eta d\zeta \approx \rho_0 \sum_{p=1}^{n_p} w_p \varphi^i(\xi_p,\eta_p,\zeta_p)\varphi^j(\xi_p,\eta_p,\zeta_p)J(\xi_p,\eta_p,\zeta_p;X_{mi})$$
$$(m=1,2,3, i=1,..,15) \tag{6}$$

Where $n_p$ stand for number of integration (Gauss) points, $w_p$ denotes weights and $\xi_p,\eta_p,\zeta_p$ are coordinates of Gauss points. Eighteen-points scheme [17] pp.80 is commonly used (e.g.[2]). For later convenience, let's consider the next definitions

$$\bar{J}_p(X_{mi}) = J(\xi_p, \eta_p, \zeta_p; X_{mi}) \ , \ \bar{M}_p^{ij} = w_p \varphi^i(\xi_p, \eta_p, \zeta_p) \varphi^j(\xi_p, \eta_p, \zeta_p) \tag{7}$$
$$(m = 1,2,3, i,j = 1,..,15, p = 1,..,n_p)$$

Such that

$$n_p = 18 \ , \ M^{ij} \approx \rho_0(\bar{J}_1 \bar{M}_1^{ij} + ... + \bar{J}_9 \bar{M}_9^{ij} + ... + \rho_0 \bar{J}_{18} \bar{M}_{18}^{ij}) \ , \ (i,j = 1,..,15) \tag{8}$$

Coefficient matrices components $\bar{M}_p^{ij}$ $(i,j=1,..,15, p=1,..,n_p)$ are scalars independent of nodal location, or solely weights. Obviously, the fewer integration points one uses, the less expensive his scheme is. How can one reduce the number of integration points in (8) and still increase the accuracy?

According to Hanukah [14-16], metric interpolation of the form $J \approx \hat{\varphi}^k \hat{J}_k$ $(k=1,..,\hat{n}_p)$, where $\hat{J}_k$ stand for metric evaluation at points $\hat{J}_k = J(\hat{\xi}_k, \hat{\eta}_k, \hat{\zeta}_k; X_{mi})$ is developed. Then, consistent mass matrix (5) rewritten as $M^{ij} \approx \rho_0 \hat{M}_k^{ij} \hat{J}_k$, $(k=1,..,\hat{n}_p)$, where coefficient scalar matrices are given by $\hat{M}_k^{ij} = \int_V \varphi^i \varphi^j \hat{\varphi}^k d\xi d\eta d\zeta$. Symmetric coefficient matrices $\hat{M}_k^{ij}$ are easily computed by MAPLE™, as they are merely integration of polynomials. In terms of (8), the new scheme is similar to $\hat{n}_p$ point rule. Herein, we develop a simple quadratic metric interpolation using which together with analytical integration result in 10-point scheme which significantly over-performs in accuracy currently used 18-point numerical integration rule.

### 3. Metric interpolation.

Jacobian matrix (3) is quadratic with respect to coordinates and can be represented by

$$J_{mn} = J_{mn}^0 + \xi J_{mn}^1 + \eta J_{mn}^2 + \zeta J_{mn}^3 + \xi\eta J_{mn}^4 + \xi\zeta J_{mn}^5 + \eta\zeta J_{mn}^6 + \xi^2 J_{mn}^7 + \eta^2 J_{mn}^8 + \zeta^2 J_{mn}^9 \tag{9}$$
$$J_{mn}^k = J_{mn}^k(X_{ri}) \ , \ (m,n,r = 1,2,3, k=0,..,9, i=1,..,15)$$

Determinant of the above Jacobian matrix is symbolically computed as a function of coordinates and nodal locations. With the help of Taylor's multivariable expansion about the origin, polynomial nature of the metric function is demonstrated

$$J = J_{order1} + J_{order2} + J_{order3} + J_{order4} + J_{order5} + J_{order6} \tag{10}$$

$$J_{order1} = J_0 + \xi J_1 + \eta J_2 + \zeta J_3$$
$$J_{order2} = \xi^2 J_4 + \xi\eta J_5 + \eta^2 J_6 + \xi\zeta J_7 + \eta\zeta J_8 + \zeta^2 J_9$$
$$J_{order3} = \xi^3 J_{10} + \xi^2\eta J_{11} + \xi\eta^2 J_{12} + \eta^3 J_{13} + \xi^2\zeta J_{14} + \xi\eta\zeta J_{15} + \eta^2\zeta J_{16} + \xi\zeta^2 J_{17} + \eta\zeta^2 J_{18} + \tag{11}$$
$$\zeta^3 J_{19}$$

$$J_{order4} = \xi^4 J_{20} + \xi^3\eta J_{21} + \xi^2\eta^2 J_{22} + \xi\eta^3 J_{23} + \eta^4 J_{24} + \xi^3\zeta J_{25} + \xi^2\eta\zeta J_{26} + \xi\eta^2\zeta J_{27} + \\ \eta^3\zeta J_{28} + \xi^2\zeta^2 J_{29} + \xi\eta\zeta^2 J_{30} + \eta^2\zeta^2 J_{31} + \xi\zeta^3 J_{32} + \eta\zeta^3 J_{33} + \zeta^4 J_{34} \quad (12)$$

$$J_{order5} = \xi^4\zeta J_{35} + \xi^3\eta\zeta J_{36} + \xi^2\eta^2\zeta J_{37} + \xi\eta^3\zeta J_{38} + \eta^4\zeta J_{39} + \xi^3\zeta^2 J_{40} + \xi^2\eta\zeta^2 J_{41} + \\ \xi\eta^2\zeta^2 J_{42} + \eta^3\zeta^2 J_{43} + \xi^2\zeta^3 J_{44} + \xi\eta\zeta^3 J_{45} + \eta^2\zeta^3 J_{46} + \xi\zeta^4 J_{47} + \eta\zeta^4 J_{48} + \zeta^5 J_{49} \quad (13)$$

$$J_{order6} = \xi^3\zeta^3 J_{50} + \xi^2\eta\zeta^3 J_{51} + \xi\eta^2\zeta^3 J_{52} + \eta^3\zeta^3 J_{53} + \xi^2\zeta^4 J_{54} + \xi\eta\zeta^4 J_{55} + \eta^2\zeta^4 J_{56} \quad (14)$$

$$\mathbf{X}_0 = \mathbf{X}(\xi=0, \eta=0, \zeta=0) \; , \; J_0(X_{mi}) = J(\mathbf{X}_0) = J(0,0,0; X_{mi}) \; , \; (m=1,2,3, i=1,..,15)$$
$$J_1(X_{mi}) = \frac{\partial J(\mathbf{X}_0)}{\partial \xi} \; , \; ... \; , \; J_9(X_{mi}) = \frac{\partial^2 J(\mathbf{X}_0)}{2\partial\zeta^2} \; , \; ... \; , \; J_{56}(X_{mi}) = \frac{\partial^6 J(\mathbf{X}_0)}{48\partial^2\eta\partial^4\zeta} \quad (15)$$

Apparently, metric (10) is rather lengthy function of nodal components and local coordinates. Partial derivatives $J_k (k=0,..,56)$ have been explicitly symbolically calculated. Importantly, exact metric representation (10) takes simple polynomial form using which, together with (5) an exact consistent mass matrix is computed and used later for error evaluation. To avoid using metric (10) with 57 precomputed terms, approximate interpolations are suggested.

The poorest (zero order) approximation for the metric is a Constant Metric CM evaluated at the centroid

$$J \approx J_{CM} = \varphi_{CM}^1 J_1^{CM} \; , \; \varphi_{CM}^1 = 1 \; , \; J_1^{CM} = \left|J_{mn}(\xi=\frac{1}{3}, \eta=\frac{1}{3}, \zeta=0)\right| \quad (16)$$

Next, Linear Metric LM interpolation is given by

$$J \approx J_{LM} = \varphi_{LM}^k J_k^{LM} \; , \; J_k^{LM} = \left|J_{mn}(p_k)\right| \; , \; (k=1,..,4, m,n=1,2,3) \quad (17)$$

$$\varphi_{LM}^1 = \frac{83}{12} - 10\xi - 10\eta - 10\zeta \; , \; \varphi_{LM}^2 = -\frac{37}{12} + 10\xi \; , \; \varphi_{LM}^3 = -\frac{37}{12} + 10\eta \; , \; \varphi_{LM}^4 = \frac{1}{4} + 10\zeta \quad (18)$$

$$p_1 : (\xi = \frac{37}{120}, \eta = \frac{37}{120}, \zeta = -\frac{1}{40}) \; , \; p_2 : (\xi = \frac{49}{120}, \eta = \frac{37}{120}, \zeta = -\frac{1}{40})$$
$$p_3 : (\xi = \frac{37}{120}, \eta = \frac{49}{120}, \zeta = -\frac{1}{40}) \; , \; p_4 : (\xi = \frac{37}{120}, \eta = \frac{37}{120}, \zeta = \frac{3}{40}) \quad (19)$$

The above uses metric evaluation at four points $p_k (k=1,..,4)$. The proposed Quadratic Metric QM interpolation is given by

$$J \approx J_{QM} = \varphi_{QM}^k J_k^{QM} \; , \; J_k^{QM} = \left|J_{mn}(p_k)\right| \; , \; (k=1,..,10, m,n=1,2,3) \quad (20)$$

Where the first four evaluation points $p_k$ are given by (19), additional six points follows

$$p_5 : (\xi = \frac{43}{120}, \eta = \frac{37}{120}, \zeta = -\frac{1}{40}) \quad , \quad p_6 : (\xi = \frac{43}{120}, \eta = \frac{43}{120}, \zeta = -\frac{1}{40})$$

$$p_7 : (\xi = \frac{37}{120}, \eta = \frac{43}{120}, \zeta = -\frac{1}{40}) \quad , \quad p_8 : (\xi = \frac{37}{120}, \eta = \frac{37}{120}, \zeta = \frac{1}{40}) \quad (21)$$

$$p_9 : (\xi = \frac{43}{120}, \eta = \frac{37}{120}, \zeta = \frac{1}{40}) \quad , \quad p_{10} : (\xi = \frac{37}{120}, \eta = \frac{43}{120}, \zeta = \frac{1}{40})$$

Quadratic functions $\varphi_{QM}^k$ (k = 1,..,10) used in (20) are given by

$$\varphi_{QM}^1 = \frac{6391}{72} - \frac{800}{3}\xi - \frac{800}{3}\eta - \frac{800}{3}\zeta + 400\xi\eta + 400\xi\zeta + 400\eta\zeta + 200\xi^2 + 200\eta^2 + 200\zeta^2$$

$$\varphi_{QM}^2 = \frac{1591}{72} - \frac{400}{3}\xi + 200\xi^2 \quad , \quad \varphi_{QM}^3 = \frac{1591}{72} - \frac{400}{3}\eta + 200\eta^2 \quad , \quad \varphi_{QM}^4 = -\frac{1}{8} + 200\zeta^2$$

$$\varphi_{QM}^5 = -\frac{3071}{36} + 400\xi + \frac{370}{3}\eta + \frac{370}{3}\zeta - 400\xi\eta - 400\xi\zeta - 400\xi^2$$

$$\varphi_{QM}^6 = \frac{1369}{36} - \frac{370}{3}\xi - \frac{370}{3}\eta + 400\xi\eta \quad (22)$$

$$\varphi_{QM}^7 = -\frac{3071}{36} + \frac{370}{3}\xi + 400\eta + \frac{370}{3}\zeta - 400\xi\eta - 400\eta\zeta - 400\eta^2$$

$$\varphi_{QM}^8 = \frac{83}{12} - 10\xi - 10\eta + \frac{800}{3}\zeta - 400\xi\zeta - 400\eta\zeta - 400\zeta^2$$

$$\varphi_{QM}^9 = -\frac{37}{12} + 10\xi - \frac{370}{3}\zeta + 400\xi\zeta \quad , \quad \varphi_{QM}^9 = -\frac{37}{12} + 10\eta - \frac{370}{3}\zeta + 400\eta\zeta$$

Importantly, CM metric (16), LM metric (17) and QM metric (20) are not unique, i.e. different sets of evaluation points and shape functions can be proposed. Moreover, we offer no rigorous optimality proof for the suggested interpolations; however, these approximations are systematic and yield satisfying accuracy. Wedge's centroid was chosen as the most representing point for CM, while points $p_k$ (k = 1,..,10) form straight-sided tetrahedron with edge length 1/10 and centroid equal to wedge's centroid. First four points are head corners, while additional six are middle edges, in fact, even higher orders interpolations can be formulated by adding evaluation point on this tetrahedron and defining shape functions accordingly, this is well known from p-version finite element method. In relation to our non-presented findings, the smaller the tetrahedron's edge length is the closer resulting interpolation to Taylor's multivariable expansion of the metric about the centroid.

## 4. Results.

CM approximation (16) together with analytical integration (5) and shape functions definition (1) result in

$$M_{CM}^{ij} = \rho_0 J_1^{CM} M_1^{ij} , \quad M_1^{ij} = \int_0^{+1} \int_0^{+1} \int_0^{1-\eta} \varphi^i \varphi^j \varphi_{CM}^1 d\xi d\eta d\zeta$$

$$M_1^{ij} = \frac{1}{1080} \begin{bmatrix} 24 & . & . & . & -30 \\ . & . & & & . \\ . & & . & & . \\ . & & & . & . \\ -30 & . & . & . & 96 \end{bmatrix} , \quad (i,j=1,..,15) \tag{23}$$

Details of symmetric coefficient matrix $M_1^{ij}$ are omitted for brevity. Please contact the corresponding author to receive any coefficient matrix necessary in either FORTRAN, text, or a MAPLE format. Following numerical integration representation (8), CM mass matrix (23) is equivalent to one point scheme. LM metric approximation (17) together with (5)(1)(18) lead to

$$M_{LM}^{ij} = \rho_0 M_k^{ij} J_k^{LM} , \quad (k=1,..,4)$$

$$M_k^{ij} = \int_0^{+1} \int_0^{+1} \int_0^{1-\eta} \varphi^i \varphi^j \varphi_{LM}^k d\xi d\eta d\zeta , \quad (i,j=1,..,15) \tag{24}$$

Once more, in terms of (8), LM mass matrix is similar (in terms of computations) to four point rules. QM approximation (20) together with analytical integration (5), shape functions (1)(22) yield

$$M_{QM}^{ij} = \rho_0 M_k^{ij} J_k^{QM} , \quad (k=1,..,10)$$

$$M_k^{ij} = \int_0^{+1} \int_0^{+1} \int_0^{1-\eta} \varphi^i \varphi^j \varphi_{QM}^k d\xi d\eta d\zeta , \quad (i,j=1,..,15) \tag{25}$$

Following (8), QM matrix (25) is computationally analogous to ten-points numerical.

5. **Preliminary numerical study.**
Consider the next element

$X_{11}=0 \quad X_{21}=0 \quad X_{31}=-1 \quad X_{12}=1 \quad X_{22}=0 \quad X_{32}=-1 \quad X_{13}=0 \quad X_{23}=1 \quad X_{33}=-1$
$X_{14}=0 \quad X_{24}=0 \quad X_{34}=1 \quad X_{15}=1 \quad X_{25}=0 \quad X_{35}=1 \quad X_{16}=0 \quad X_{26}=1 \quad X_{36}=1$
$X_{17}=\frac{1}{2} \quad X_{27}=0 \quad X_{37}=-1 \quad X_{18}=\frac{1}{2} \quad X_{28}=\frac{1}{2} \quad X_{38}=-1 \quad X_{19}=0 \quad X_{29}=\frac{1}{2} \quad X_{39}=-1$ (26)
$X_{110}=\frac{1}{2} \quad X_{210}=0 \quad X_{310}=1 \quad X_{111}=\frac{1}{2} \quad X_{211}=\frac{1}{2} \quad X_{311}=1 \quad X_{112}=0 \quad X_{212}=\frac{1}{2} \quad X_{312}=1$
$X_{113}=0 \quad X_{213}=0 \quad X_{313}=0 \quad X_{114}=1 \quad X_{214}=0 \quad X_{314}=0 \quad X_{115}=0 \quad X_{215}=1 \quad X_{315}=0$

The above is a fine mesh element (parent element) with constant metric. We introduce the coarseness measure $\delta$ by replacing several components such that for $\delta = 0$ result in fine element (26), while for $\delta > 0$ yield coarse mesh element. Consider the next three element families

$$X_{34}=1+\delta, \ X_{15}=1+\delta, \ X_{26}=1+\delta \qquad (27)$$

$$X_{11}=0+\delta, \ X_{22}=0+\delta, \ X_{32}=-1+\delta \qquad (28)$$

$$X_{34}=1-\delta, \ X_{35}=1-\delta, \ X_{310}=1-\delta \qquad (29)$$

The metric of (27) is a six order expression $J = 1 - 1.5\delta + 0.5\delta^2 + ... + 1.5\delta^3\eta^2\zeta^3$, (28) produce fourth order metric $J = 1 - \delta - 1.5\delta^2 + ... - \delta^2\zeta^4$, and (29) leads to third order metric $J = 1 - 0.5\delta + ... - \delta^2\zeta^3$.

Coarseness $\delta$ has been gradually increases. For each value of $\delta$, CM (23), LM (24) and QM (25) mass matrices have been calculated, absolute error for each components is computed with respect to the exact values, and an averaged (among all the components) absolute error values are recorded at Figure 2.

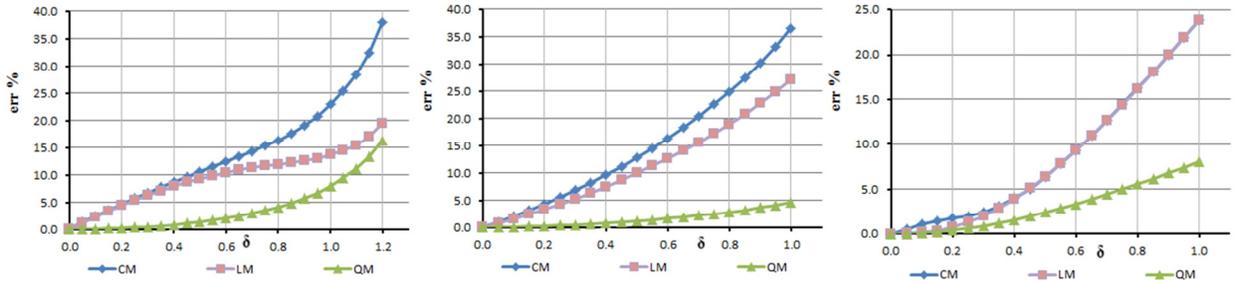

Figure 2: An averaged absolute error of CM, LM and QM mass matrices are reported as a function of the mesh coarseness $\delta$. Left graph records the first element (27) results; middle graph corresponds to (28) while the right graph contains third element (29) details.

According to expectations, QM is the most accurate, LM over-performs CM. It is recalled that CM is one point scheme, LM equivalent to four-point numerical integration and QM is computationally similar to ten point integration rule.

An averaged absolute error of QM is compared to 18-point numerical integration Figure 3.

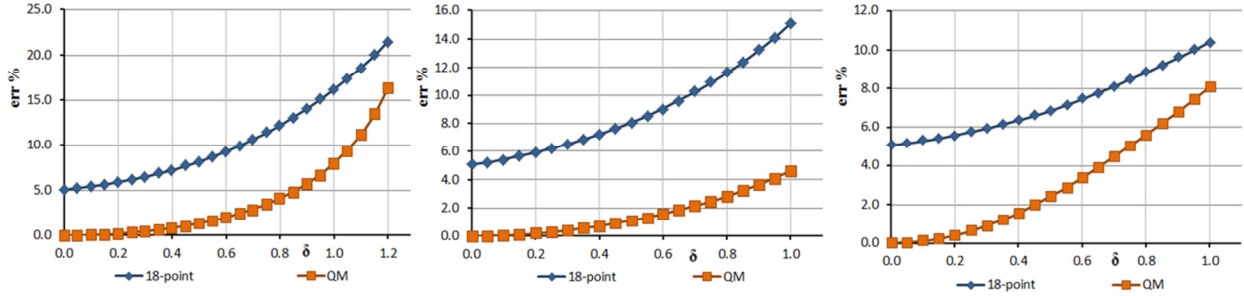

Figure 3: An averaged absolute error of QM and 18-point scheme. Left graph for (27), then the second element (28) in the middle, on the right is the third element family (29).

Clearly, for all the coarseness range $\delta$, QM over-performs 18-point scheme. Next we examine the maximum absolute error (among all the components) of QM and 18-point numerical integration Figure 4.

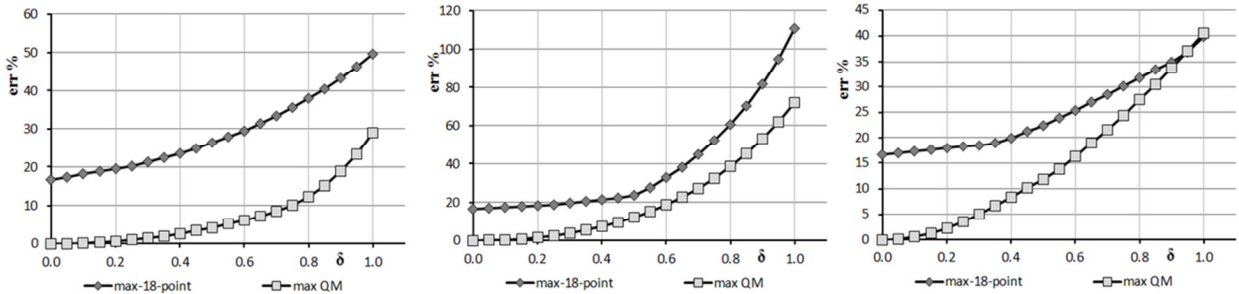

Figure 4: Maximum absolute error for QM vs 18-point rule. Left graph corresponds to (27), then graph for an element family (28) and the right graph stands for (29).

In terms of maximum absolute error, our QM is superior to 18-point scheme. The conducted *preliminary* accuracy study, though incomplete, still clearly demonstrates potential benefits of the suggested QM formulation.

## 6. Conclusions

New, ten-point integration rule for mass matrix of fifteen-node wedge element (3D solid) have been developed. To this end, special second order metric interpolation requiring ten evaluation points is proposed. Later, analytical integration for coefficient matrices (weights) led to ten-point integration rule. Preliminary numerical study including fine and gradually changing coarse mesh has been conducted. For each coarseness $\delta$ value, maximum absolute error (among all the components) and an averaged absolute error has been calculated and presented in graphical form. Certainly, according to our findings, our ten point rule superior in accuracy both for maximum and for an averaged absolute error.